\documentclass[11pt,leqno,a4paper]{article}

\usepackage{amsmath, amssymb, amsthm, amsfonts}
\usepackage{mathrsfs}
\usepackage{color}
\theoremstyle{plain}
\newtheorem{thm}{Theorem}

\newtheorem{lem}{Lemma}

\theoremstyle{remark}
\newtheorem{rem}{Remark}

\def\K{\mathop{\mbox{\bf\Large K}}}

\numberwithin{equation}{section}

\begin{document}

\title{Multiple-Correction and Continued Fraction Approximation}
\author{Xiaodong Cao}
\date{}

\maketitle

\footnote[0]{2010 Mathematics Subject Classification
Primary 11Y60 41A25 34E05 26D15}

\footnote[0]{This work is supported by
the National Natural Science Foundation of China (Grant No.11171344) and the Natural
Science Foundation of Beijing (Grant No.1112010).}

\begin{abstract}
The main aim of this paper is to further develop a multiple-correction
method formulated in a previous work~\cite{CXY}. As its applications, we find a kind of hybrid-type finite
continued fraction approximations in two cases of Landau constants and Lebesgue constants. In addition, we refine the previous results of Lu~\cite{Lu} and Xu and You~\cite{XY} for the Euler-Mascheroni constant.

\end{abstract}

\section{Introduction}

The constants of Landau and Lebesgue are defined for all integers $n\ge 0$, respectively, by
\begin{align}
G(n)=\sum_{k=0}^{n}\frac{1}{16^k}\binom {2k}{k}^2\quad  \mbox{and}\quad
L_n=\frac{1}{2\pi}\int_{-\pi}^{\pi}\left|\frac{\sin\left((n+\frac 12)t\right)}{\sin(\frac t2)}\right|dt.
\end{align}
The constants $G(n)$ are important in complex analysis.
In 1913, Landau~\cite{La} proved that if $f(z)=\sum_{k=0}^{\infty}a_kz^k$ is an analytic function in the unit disc satisfying
$|f(z)|<1$ for $|z|<1$, then $\left|\sum_{k=0}^{n}a_k\right|\le G(n)$, and that this bound is optimal. Furthermore, Landau~\cite{La} showed that
\begin{align}
G(n)\sim \frac{1}{\pi}\ln n, (n\rightarrow \infty).
\label{Landau-C-1}
\end{align}
Let $\gamma$ denote the Euler-Mascheroni constant. In 1930, Watson~\cite{Wat} obtained a more precise asymptotic
formula than~\eqref{Landau-C-1}
\begin{align}
G(n)\sim \frac{1}{\pi}\ln (n+1)+c_0-
\frac{1}{4\pi (n+1)}+O\left(\frac{1}{n^2}\right),
(n\rightarrow \infty),\label{Landau-C-2}
\end{align}
where
\begin{align}
c_0=\frac{1}{\pi}(\gamma+4\ln 2)=1.0662758532089143543\cdots.
\label{c0-definition}
\end{align}
In fact, the work of Watson opened up a novel insight into the asymptotic behavior of the Landau sequences $(G(n))_{n\ge0}$. Inspired by~\eqref{Landau-C-2}, many authors investigated the
upper and lower bounds of $G(n)$. We list some main results as follows:
\begin{eqnarray}
~~&&\frac{1}{\pi}\ln (n+1)+1\le G(n)
<\frac{1}{\pi}\ln (n+1)+c_0\quad (n\ge 0),\quad\mbox{(Brutman~\cite{Br},1982)}\\
~~&&\frac{1}{\pi}\ln \left(n+\frac 34\right)+c_0<G(n)
\le\frac{1}{\pi}\ln\left(n+\frac 34\right)+1.0976\quad (n\ge 0),\quad\mbox{(Falaleev~\cite{Fal}, 1991)}\label{Landau-Falaleev}\\
~~&&\frac{1}{\pi}\ln \left(n+\frac 34\right)+c_0<G(n)
<\frac{1}{\pi}\ln\left(n+\frac 34+\frac{11}{192n}\right)+c_0\quad (n\ge 1),\quad\mbox{(Mortici~\cite{Mor4}. 2011)}\label{Landau-Mortici}
\end{eqnarray}
Recently, Chen~\cite{Ch2} found the following better
approximation for $G(n)$: as $n\rightarrow\infty$,
\begin{align}
G(n)
=&c_0+\frac{1}{\pi}\ln\left(n+\frac 34+\frac{11}{192(n+\frac 34)}-\frac{2009}{184320(n+\frac 34)^3}
+\frac{2599153}{371589(n+\frac 34)^5}\right)\\
&+O\left(\frac{1}{(n+\frac 34)^8}\right),\nonumber
\end{align}
and the better upper bound:
\begin{align}
G(n)
<c_0+\frac{1}{\pi}\ln\left(n+\frac 34+\frac{11}{192(n+\frac 34)}\right), (n\ge 0).\label{chen-1}
\end{align}
More recently, Cao, Xu and You~\cite{CXY} improved the rate of convergence to $n^{-14}$, and attained the following tight double-sides inequalities
\begin{align}
\frac{C_1}{(n+\frac 32)^6}
<G(n)-\frac{1}{\pi}\ln(n+\frac 34)-c_0-\frac{\frac{11}{192\pi}}{(n+\frac 34)^2+\frac{1541}{7040}}
<\frac{C_1}{(n+\frac 12)^6},(n\ge 0),
\end{align}
where $C_1=\frac{89684299}{18166579200\pi}$.

Another direction for developing the approximation to $G(n)$ was
initiated by Cvijovi\'c and Klinowski~\cite{CK}, who established the following estimates
of $G(n)$ in terms of the Psi(or Digamma) function $\psi(z):=\frac{\Gamma'(z)}{\Gamma(z)}$:
\begin{align}
&\frac{1}{\pi}\psi\left(n+\frac 54\right)+c_0<G(n)
<\frac{1}{\pi}\psi\left(n+\frac 54\right)+1.0725,\quad(n\ge 0),\\
&\frac{1}{\pi}\psi\left(n+\frac 32\right)+0.9883<G(n)
<\frac{1}{\pi}\psi\left(n+\frac 32\right)+c_0,\quad(n\ge 0).
\end{align}
Since then, many authors have made significant contributions to sharper
the inequalities and the asymptotic expansions for $G(n)$, see e.g.
Alzer~\cite{Al}, Chen~\cite{Ch1}, Cvijovi\'c and Srivastava~\cite{CS}, Granath~\cite{Gra}, Mortici~\cite{Mor4}, Nemes~\cite{Ne1,Ne2}, Popa~\cite{Po}, Popa and Secelean~\cite{PS}, Zhao~\cite{Zhao},  Gavrea and M. Ivan~\cite{GI}, Chen and Choi~\cite{CC1,Ch2,CC2}, etc. To the best knowledge of the authors, the latest lower and upper bounds of $G(n)$ along this research direction are due to Chen and Choi~\cite{CC2}.

In 1906, Lebesgue~\cite{Le} showed that if a function $f$ is integrable on the interval $[-\pi,\pi]$
 and $S_n(f;x)$ is the $n$-th partial sum of the Fourier series of $f$, then, we have
 \begin{align*}
 &a_k=\frac{1}{\pi}\int_{-\pi}^{\pi}f(t)\cos(kt)dt\quad\mbox{($k\ge 0$) and}
 \quad  b_k=\frac{1}{\pi}\int_{-\pi}^{\pi}f(t)\sin(kt)dt,\\
 &S_n(f;x)=\frac{a_0}{2}+\sum_{k=1}^{n}\left(a_k \cos(kx)
 +b_k \sin(kx)
 \right),
 \end{align*}
 where the sum for $n=0$ is usually stipulated to be zero. If $|f(x)|\le 1$ for all
 $x\in [-\pi,\pi]$, then
 \begin{align}
 \max_{x\in [-\pi,\pi]}|S_n(f;x)|\le L_n\quad (n\ge 0).
 \label{Lebesgue-definition-1}
 \end{align}
It is noted that $L_n$ is the smallest possible constant for which
the inequality~\eqref{Lebesgue-definition-1} holds for all
integrable functions $f$ on $[-\pi,\pi]$.

The Lebesgue constants play an important role in the theory of Fourier
series. Therefore, they have attracted much attention of several well-known mathematicians such as  Fej\'er~\cite{Fe}, Gronwall~\cite{Gro}, Hardy~\cite{Ha}, Szeg\"o~\cite{Sz}, Watson~\cite{Wat}, who established some remarkable properties of these
numbers including monotonicity theorems, and various series and integral representations for $L_n$. Watson~\cite{Wat} showed
\begin{align}
L_{n/2}=\frac{4}{\pi^2}\ln (n+1)+c_1+O\left(\frac{1}{n^2}\right)\quad (n\rightarrow\infty),
\label{Lebesgue-watson}
\end{align}
where
\begin{align}
c_1=\frac{8}{\pi^2}\sum_{k=1}^{\infty}
\frac{\ln k}{4k^2-1}+\frac{4}{\pi^2}(\gamma+2\ln 2)=0.98943127383114695174\cdots.\label{c1-definition}
\end{align}
Since then, many authors have made important contributions to this research topic, see e.g. Galkin~\cite{Gal}, Wong~\cite{Wo}, Alzer~\cite{Al}, Zhao~\cite{Zhao}, Chen and Choi~\cite{CM}, etc.
Let
\begin{align}
u_n=&L_{n/2}-\left(c_1+\frac{4}{\pi^2}\psi\left(
n+\frac 32+\frac{\frac 18-\frac{\pi^2}{72}}{n+1}\right)\right),\\
v_n=&L_{n/2}-\left(c_1+\frac{4}{\pi^2}\ln\left(
n+1+\frac{a}{n+1}+\frac{c}{(n+1)^3}
\right)\right),
\end{align}
where
\begin{align*}
a=\frac{1}{6}-\frac{\pi^2}{72}\quad \mbox{and}\quad
c=-\frac{37}{360}+\frac{\pi^2}{135}+\frac{67\pi^4}{259200}.
\end{align*}
Recently, Chen and Choi~\cite{CC1} obtained
\begin{align}
\lim_{n\rightarrow\infty}n^4u_n=&\frac{-23625+
1770\pi^2+67\pi^4}{64800\pi^2},\\
\lim_{n\rightarrow\infty}n^6v_n=&-\frac{-188637120+
15135120\pi^2+308196\pi^4+7537\pi^6}{97977600\pi^2}.
\end{align}

\noindent {\bf Notation.} Throughout the paper, the notation
$P_k(x)$(or $Q_k(x)$) as usual denotes a polynomial of degree $k$
in terms of $x$. The notation $\Psi(k;x)$
means a polynomial of degree $k$ in terms of $x$ with all of its non-zero coefficients being positive, which may be different at each occurrence.
Notation $\Phi(k;x)$ denotes a polynomial of degree $k$ in terms of $x$ with
the leading coefficient being equal to one, which may be different at different subsections.

This paper is organized as follows. In Section 2, we prepare
some  preliminary lemmas. In Section 3, we explain how to find
a finite continued fraction approximation by using the
\emph{multiple-correction method}. In Section 4 and Section 5,
we discuss the constants of Landau and Lebesgue, respectively.
In the last section, we consider to refine the works of Lu~\cite{Lu}
and Xu and You~\cite{XY} for the Euler-Mascheroni constant.

\section{Some Preliminary Lemmas}

The following lemma gives a method for measuring the rate of convergence, whose proof can be found in~\cite{Mor1,Mor2}.
\begin{lem}
If the sequence $(x_n)_{n\in \mathbb{N}}$ is convergent to zero and there exists the limit
\begin{align}
\lim_{n\rightarrow +\infty}n^s(x_n-x_{n+1})=l\in [-\infty,+\infty]
\label{LEM1-1}
\end{align}
with $s>1$, then
\begin{align}
\lim_{n\rightarrow +\infty}n^{s-1}x_n=\frac{l}{s-1}.\label{LEM1-2}
\end{align}
\end{lem}
In the study of Landau constants, we need to apply a so-called
Brouncker's continued fraction formula.
\begin{lem}
For all integer $n\ge 0$, we have
\begin{align}
q(n):=\left(\frac{\Gamma(n+\frac 12)}{\Gamma(n+1)}\right)^2=
\frac{4}{1+4n+\frac{1^2}{2+8n+\frac{3^2}{2+8n
+\frac{5^2}{2+8n+\ddots}}}}.\label{LEM2-1}
\end{align}
\end{lem}
In 1654 when Brouncker and Wallis collaborated on the problem of squaring the circle, Lord William Brouncker found this remarkable fraction formula. Formula~\eqref{LEM2-1} was not published by Brouncker himself, and first
appeared in~\cite{Wal}. For a general $n$, it actually follows from Entry 25 in Chapter 12 in Ramanujan's notebook~\cite{Ber}, which gives a more
general continued fraction formula for quotients of gamma functions.

Writing continued fractions in the way of~\eqref{LEM2-1} takes a lot of space and thus, we use the following shorthand notation
\begin{align}
q(n)=\frac{4}{1+4n+}\frac{1^2}{2+8n+}\frac{3^2}{2+8n+}
\frac{5^2}{2+8n+}\cdots
=\frac{4}{1+4n+}\K_{k=0}^{\infty}
\frac{(2k+1)^2}{2+8n},
\end{align}
and its $k$-th approximation $q_k(n)$ is defined by
\begin{align}
q_1(n)&=\frac{4}{1+4n},\\
q_k(n)&=\frac{4}{1+4n+}\frac{1^2}{2+8n+}\frac{3^2}{2+8n+}
\frac{(2k-3)^2}{2+8n}
=\frac{4}{1+4n+}\K_{j=0}^{k-1}
\frac{(2j-1)^2}{2+8n}
,\quad (k\ge 2).\label{qk-definition}
\end{align}
\begin{lem}
Let $c_1$ be defined by~\eqref{c1-definition}. Then,
for $n\in \mathbb{N}_0$ and $N\in \mathbb{N}$, we have
\begin{align}
&\frac{4}{\pi^2}\ln (n+1)+c_1+
\sum_{j=1}^{2N}
\frac{a_j}
{(n+1)^{2j}}\\
&<L_{n/2}
<\frac{4}{\pi^2}\ln (n+1)+c_1+
\sum_{j=1}^{2N+1}
\frac{a_j}{(n+1)^{2j}},\nonumber
\end{align}
where
\begin{align*}
a_j:=\frac{8}{\pi^2}\frac{B_{2j}}{2j}\left(2^{2j-1}-1\right)
\left(1+\sum_{k=1}^{j}\frac{(-1)^k}{(2k)!}B_{2k}\pi^{2k}\right),
\end{align*}
and the Bornoulli numbers $B_k$
is defined by
\begin{align*}
\frac{z}{e^z-1}=\sum_{k=0}^{\infty}B_k\frac{z^k}{k!}\quad (|z|<2\pi).
\end{align*}
\end{lem}
\proof  This is Theorem 3.1 of Chen and Choi~\cite{CC2}, and also see~(3.8) in Chen and Choi~\cite{CC1}.\hfill$\Box$

In the proof of our inequalities for the constants of Landau and
Lebesgue, we also need to use the following simple inequality,
which is a consequence of Hermite-Hadamard inequality.
\begin{lem} Let $f$ be twice derivable, with $f''$ continuous.  If $f''(x)>0$, then
\begin{align}
\int_{a}^{a+1}f(x) dx > f(a+1/2).\label{LEM3}
\end{align}
\end{lem}

\section{The multiple-correction method}

First, let us briefly review a so-called \emph{multiple-correction method} presented in our previous paper~\cite{CXY}.
Let $(v(n))_{n\ge 1}$ be a sequence to be approximated. Throughout
the paper, we always assume that the following three
conditions hold.

\noindent {{\textbf{Condition} (i)}.} The initial-correction function $\eta_0(n)$ satisfies
\begin{align*}
&\lim_{n\rightarrow\infty}\left(v(n)-\eta_0(n)\right)=0,\\
&\lim_{n\rightarrow\infty}n^{l_0}\left(v(n)-v(n+1)-\eta_0(n)
+\eta_0(n+1)
\right)=C_0\neq 0,
\end{align*}
for some positive integer $l\ge 2$.

\noindent {{\textbf{Condition} (ii)}.} The $k$-th correction function $\eta_k(n)$ has the form of $-\frac{C_{k-1}}{\Phi_k(l_{k-1};n)}$, where
\begin{align*}
\lim_{n\rightarrow\infty}n^{l_{k-1}}\left(v(n)-v(n+1)-
\sum_{j=0}^{k-1}\left(\eta_j(n)-\eta_j(n+1)\right)
\right)=C_{k-1}\neq 0,
\end{align*}

\noindent {{\textbf{Condition} (iii)}.} The difference  $\left(v(1/x)-v(1/x+1)-\eta_0(1/x)
+\eta_0(1/x+1)
\right)$ is an analytic function in a neighborhood of point $x=0$.

Actually, the \emph{multiple-correction method} is a recursive algorithm. If the assertion
\begin{align*}
\lim_{n\rightarrow\infty}n^{l_{k-1}}\left(v(n)-v(n+1)-
\sum_{j=0}^{k-1}\left(\eta_j(n)-\eta_j(n+1)\right)
\right)=C_{k-1}\neq 0
\end{align*}
is true, then, it is not difficult to observe
\begin{align*}
\lim_{n\rightarrow\infty}\left(\left(v(n)-v(n+1)-
\sum_{j=0}^{k-1}\left(\eta_j(n)-\eta_j(n+1)\right)
\right)-\frac{C_{k-1}}{n^{l_{k-1}}}\right)=0.
\end{align*}
Roughly speaking, the idea of the \emph{multiple-correction method} is to use the polynomial $\Phi_k(l_{k-1};n)$ of degree $l_{k-1}$ instead of $n^{l_{k-1}}$ for improving the convergence rate. In other words, we view $n^{l_{k-1}}$ as a special polynomial of degree $l_{k-1}$ in terms of $n$.
Here, we note that the polynomial $\Phi_k(l_{k-1};n)$ contains  $l_{k-1}$ undetermined
parameters $a_j$ $(0\le j\le l_{k-1}-1)$, and hence, in some cases we hope that we can attain ``more gains".

The initial-correction is a very important step. With this, we hope to further develop the above method starting from the second-correction. To find
the proper structure of finite continued fraction, we must try many times
by using $-\frac{C_0}{\Phi_1(l_{1};n)+\frac{b_j}{n^j}}$ instead of $-\frac{C_0}{\Phi_1(l_{1};n)}$, where $j$ is a positive integer.
To do that, we need to begin from $j=1$ and try step by step. Once we have found that the convergence rate can be improved for the first positive integer, say $j_0$, we use $\Phi(j_0; n)$ to replace $n^{j_0}$ immediately, and then, determine all the corresponding coefficients of the polynomial $\Phi(j_0; n)$. We continue this process until the desired structure of finite continued fraction is found.
It is for this reason that we call it as the \emph{multiple-correction method}.

In addition, to determine all the related coefficients, we often use an appropriate symbolic computation software, which needs a huge of computations. On the other hand, the exact expressions at each occurrence also takes a lot of space. Hence, in this paper we omit some related details for space limitation. For interesting readers, see our previous paper~\cite{CXY}.

It is a natural question whether or not \emph{multiple-correction method} can be used to accelerate convergence in some BBP-type or Ramanujan-type series, we hope to return to this topic elsewhere.

\section{The Landau Constants}
\begin{thm} Let sequences ${\rm MC}_k(n)$ be defined as follows:
\begin{align}
\mathrm{MC}_0(n): =&\frac{1}{\pi}\ln\left(n+\frac 34\right)+c_0,\\
\mathrm{MC}_1(n):=&\frac{1}{\pi}\frac{\kappa_1}{(n+\frac 34)^2+\lambda_1},
\\
\mathrm{MC}_k(n):=&\frac{1}{\pi}\frac{\kappa_1}{(n+\frac 34)^2+\lambda_1+}
\K_{j=2}^{k}
\frac{\kappa_j}{(n+\frac 34)^2+\lambda_j},
(k\ge 2),\label{LMC-definition}
\end{align}
where $c_0$ is determined by~\eqref{c0-definition} and
\begin{align*}
\kappa_1=&\frac{11}{192},&\lambda_1=\frac{1541}{7040},\\
\kappa_2=&-\frac{89684299}{1040793600},&\lambda_2=\frac{815593360691}
{631377464960},\\
\kappa_3=&-\frac{791896453750695892475}{691850212268234428416},
&\lambda_3=\frac{79124827964452580408836456738931}
{23635681749960244849264556808320}.
\end{align*}
If we let the $k$-th correction error term $E_k(n)$ be denoted by
\begin{align}
E_k(n):=G(n)-\mathrm{MC}_0(n)-\mathrm{MC}_k(n),\label{LE-definition}
\end{align}
then, for all positive integer $k$, we have
\begin{align}
&\lim_{n\rightarrow\infty}n^{4k+3}\left(E_k(n)-E_k(n+1)\right)
=(4k+2)C_k,\\
&\lim_{n\rightarrow\infty}n^{4k+2}E_k(n)=C_k,\label{LMC-result}
\end{align}
where
 \begin{eqnarray*}
C_1&=&\frac{89684299}{18166579200\pi},\\
C_2&=&\frac{31675858150027835699}{5605686531912433139712\pi},\\
C_3&=&\frac{9662255454831353335643376823083291821}
{310776980771128296411407710029663436800\pi}.
 \end{eqnarray*}
\end{thm}

\begin{rem}
Theorem 1 tells us that it could be possible for us to find a simpler asymptotic expansion than Theorem 2.1 of Chen and Choi~\cite{CC1} for the Landau constants.
\end{rem}

\proof Let us consider the initial-correction.

\noindent{\bf (Step 1) The initial-correction.}
Motivated by inequalities \eqref{Landau-Falaleev} and \eqref{Landau-Mortici},
we choose $\mathrm{MC}_0(n)=\frac{1}{\pi}\ln(n+\frac 34)+c_0$, and define
\begin{align}
E_0(n)=G(n)-\mathrm{MC}_0(n)=G(n)-\frac{1}{\pi}\ln(n+\frac 34)-c_0.
\label{L-E0-definition}
\end{align}
Then, it follows
immediately from \eqref{L-E0-definition}
\begin{align}
E_0(n)-E_0(n+1)=G(n)-G(n+1)-\frac{1}{\pi}\ln(n+\frac 34)+
\frac{1}{\pi}\ln(n+\frac 74).\label{L-E0-difference}
\end{align}
Now, by using the duplication formula (Legendre, 1809)
\begin{align}
2^{2z-1}\Gamma(z)\Gamma(z+\frac 12)=\sqrt\pi\Gamma(2z),
\label{duplication formula}
\end{align}
one can prove
\begin{align}
G(n)-G(n-1)=\frac{(\Gamma(2n+1))^2}{16^2(\Gamma(n+1))^4}
=\left(\frac{(2n)!}{4^n(n!)^2}\right)^2=\frac{1}{\pi}q(n),
\label{G-difference-def}
\end{align}
where $q(n)$ is defined by~\eqref{LEM2-1}. Also see Page 739 in Granath~\cite{Gra} or Page 306 in Chen~\cite{Ch2}.
Combining~\eqref{L-E0-difference} with \eqref{G-difference-def} yields
\begin{align}
E_0(n)-E_0(n+1)=-\frac{1}{\pi}q(n+1)-\frac{1}{\pi}\ln(n+\frac 34)+
\frac{1}{\pi}\ln(n+\frac 74).\label{L-E0-difference-1}
\end{align}
On one hand, by utilizing Lemma 2 and~\eqref{qk-definition}, we can obtain that for all positive integer $j$,
\begin{align}
q_2(n)<q_4(n)<\cdots<q_{2j}(n)<q(n)
<q_{2j+1}(n)<\cdots<q_3(n)<q_1(n).\label{Landau-1}
\end{align}
On the other hand,  by using \emph{Mathematica} software, we can attain
\begin{align}
q_9(n)-q_8(n)=O\left(\frac{1}{n^{17}}\right).\label{Landau-1-2}
\end{align}
Now, combining \eqref{Landau-1} and \eqref{Landau-1-2} gives us
\begin{align}
q(n+1)=q_{8}(n+1)+O\left(\frac{1}{n^{16}}\right).
\label{Landau-1-3}
\end{align}
Again, by making use of the \emph{Mathematica} software, we can
expand $q_{8}(n+1)$ into a power series in terms of $n^{-1}$ so that
\begin{align}
q(n+1)=&q_{8}(n+1)+O\left(\frac{1}{n^{16}}
\right)\label{q8-approximation}\\
=&\frac{1}{n}-\frac 54\frac{1}{n^2}+\frac{49}{32}\frac{1}{n^3}
-\frac{235}{128}\frac{1}{n^4}+\frac{4411}{2048}\frac{1}{n^5}
-\frac{20275}{8192}\frac{1}{n^6}
+\frac{183077}{65536}\frac{1}{n^7}\nonumber\\
&-\frac{815195}{262144}\frac{1}{n^8}+
\frac{28754131}{8388608}\frac{1}{n^9}
-\frac{125799895}{33554432}\frac{1}{n^{10}}
+\frac{1091975567}{268435456}\frac{1}{n^{11}}\nonumber\\
&-\frac{4702048685}{1073741824}\frac{1}{n^{12}}
+\frac{80679143663}{17179869184}\frac{1}{n^{13}}-
\frac{346250976095}{68719476736}\frac{1}{n^{14}}\nonumber\\
&+\frac{2947620308941}{549755813888}\frac{1}{n^{15}}
+O\left(\frac{1}{n^{16}}\right).\nonumber
\end{align}
In addition, it is not difficult to obtain
\begin{align}
-\ln(n+\frac 34)+
\ln(n+\frac 74)=\frac{1}{n}-\frac{5}{4}\frac{1}{n^2}
+\frac{79}{48}\frac{1}{n^3}+O\left(\frac{1}{n^{4}}\right).
\label{Landau-1-4}
\end{align}
Inserting \eqref{q8-approximation} and \eqref{Landau-1-4} into
\eqref{L-E0-difference-1} results in
\begin{align}
E_0(n)-E_0(n+1)=\frac{11}{96\pi}\frac{1}{n^3}
+O\left(\frac{1}{n^{4}}\right).\label{Landau-1-5}
\end{align}
Note that the inequalities \eqref{Landau-Mortici} implies $E_0(\infty)=0$. By Lemma 1 again, we obtain
\begin{align}
\lim_{n\rightarrow\infty}n^2E_0(n)=\frac{11}{192\pi}=C_0=\kappa_1.
\label{Landau-1-6}
\end{align}

\noindent{\bf (Step 2) The first-correction.} For simplicity, let
\begin{align}
\mathrm{MC}_1(n)=\frac{1}{\pi}\frac{\kappa_1}{\Phi_1(2;n)}
=\frac{1}{\pi}\frac{\kappa_1}{(n+\frac 34)^2+\lambda_1},\label{Landau2-1}
\end{align}
and define
\begin{align}
E_1(n):=G(n)-\mathrm{MC}_0(n)-\mathrm{MC}_1(n)
=E_0(n)-\mathrm{MC}_1(n).\label{Landau2-2}
\end{align}
Combining~\eqref{L-E0-difference-1}, \eqref{Landau-1-3} and~\eqref{Landau2-2}, we can obtain
\begin{align}
E_1(n)-E_1(n+1)=&\left(E_0(n)-E_0(n+1)\right)
-\left(\mathrm{MC}_1(n)-\mathrm{MC}_1(n+1)\right)\label{Landau2-3}\\
=&-\frac{1}{\pi}q_8(n+1)-\frac{1}{\pi}\ln(n+\frac 34)+
\frac{1}{\pi}\ln(n+\frac 74)\nonumber\\
&-\mathrm{MC}_1(n)+\mathrm{MC}_1(n+1)
+O\left(\frac{1}{n^{16}}\right).\nonumber
\end{align}
By taking advantage of formulae~\eqref{q8-approximation}
and \eqref{Landau2-1}, and~\emph{Mathematica} software, we expand $E_1(n)-E_1(n+1)$ into power series in terms of $n^{-1}$:
\begin{align}
\pi\left(E_1(n)-E_1(n+1)\right)=
&\frac{-\frac{1541}{30720}+\frac{11\lambda_1}{48}}{n^5}
+\frac{\frac{7705}{24576}-\frac{275 \lambda_1}{192}}{n^6}
\label{Landau2-4}\\
&+\frac{-\frac{183077}{65536}+\frac{275463+973280\lambda_1-59136 \lambda_1^2}{172032}}{n^7}+O\left(\frac{1}{n^{8}}\right).
\nonumber
\end{align}
By Lemma 1, the fastest sequence $(E_1(n))_{n\ge 1}$ is obtained when the first
coefficient of this power series vanish. In this case
\begin{align}
\lambda_1=\frac{1541}{7040},\label{Landau2-5}
\end{align}
and thus,
\begin{align*}
E_1(n)-E_1(n+1)=\frac{89684299}{3027763200\pi}\frac{1}{n^7}+
O\left(\frac{1}{n^{8}}\right).
\end{align*}
Now, by Lemma 1 again, we attain
\begin{align}
\lim_{n\rightarrow\infty}n^6E_1(n)=
\frac{89684299}{18166579200\pi}=C_1.
\label{Landau2-8}
\end{align}

\noindent{\bf (Step 3) The second-correction.} Let
\begin{align}
\mathrm{MC}_2(n)=\frac{1}{\pi}\frac{\kappa_1}{(n+\frac 34)^2+\lambda_1+}\frac{\kappa_2}{(n+\frac 34)^2+\lambda_2},\label{Landau3-1}
\end{align}
and define
\begin{align}
E_2(n)=G(n)-\mathrm{MC}_0(n)-\mathrm{MC}_2(n).\label{Landau3-2}
\end{align}
Following the way similar to the proof of~\eqref{Landau2-3}, we can obtain
\begin{align}
E_2(n)-E_2(n+1)=&-\frac{1}{\pi}q_8(n+1)-\frac{1}{\pi}\ln(n+\frac 34)+
\frac{1}{\pi}\ln(n+\frac 74)\label{Landau3-4}\\
&-\mathrm{MC}_2(n)+\mathrm{MC}_2(n+1)
+O\left(\frac{1}{n^{16}}\right).\nonumber
\end{align}
By utilizing~\emph{Mathematica} software, $E_2(n)-E_2(n+1)$ can be expanded into power series in terms of $n^{-1}$
\begin{align}
\pi\left(E_2(n)-E_2(n+1)\right)=&\frac{\frac{89684299}{3027763200}+
\frac{11\kappa_2}{32}}{n^7}
-\frac{\frac{89684299}{346030080}+\frac{385\kappa_2}{128}}{n^8}
\label{Landau3-5}\\
&+\frac{\frac{2961426180353}{2283798528000}+\frac{120119\kappa_2}{7680}
-\frac{11\kappa_2\lambda_2}{24}}{n^9}\nonumber\\
&+\frac{-\frac{988371602353}{203004313600}
-\frac{129357\kappa_2}{2048}
+\frac{165\kappa_2\lambda_2}{
32}}{n^{10}}+O\left(\frac{1}{n^{12}}\right)\nonumber\\
&+\frac{\frac{65353200785578639}{4287451103232000}+\frac{6308113241
\kappa_2}{28835840}-\frac{207019\kappa_2\lambda_2}{6144}+
\frac{55\kappa_2\lambda_2^2}{96}-\frac{55\kappa_2^2}{96}}
{n^{11}}.\nonumber
\end{align}
The fastest sequence $(E_2(n))_{n\ge 1}$ is obtained by enforcing the first
four coefficients of this power series to be zeros. In this case
\begin{align}
\kappa_2=-\frac{89684299}{1040793600}\quad\mbox{and}\quad
\lambda_2=\frac{815593360691}{631377464960},\label{Landau3-6}
\end{align}
and hence,
\begin{align*}
E_2(n)-E_2(n+1)=\frac{158379290750139178495}
{2802843265956216569856\pi}\frac{1}{n^{11}}+
O\left(\frac{1}{n^{12}}\right).
\end{align*}
Combining this with Lemma 1 leads to
\begin{align}
\lim_{n\rightarrow\infty}n^{10}E_2(n)=
\frac{31675858150027835699}{5605686531912433139712\pi}:=C_2.
\label{Landau3-7}
\end{align}

\noindent{\bf (Step 4) The third-correction.} Let
\begin{align}
\mathrm{MC}_3(n):=\frac{1}{\pi}\frac{\kappa_1}{(n+\frac 34)^2+\lambda_1+}
\frac{\kappa_2}{(n+\frac 34)^2+\lambda_2+}
\frac{\kappa_3}{(n+\frac 34)^2+\lambda_3},
\label{Landau4-1}
\end{align}
If we define
\begin{align}
E_3(n)=G(n)-\mathrm{MC}_0(n)-\mathrm{MC}_3(n).\label{Landau4-2}
\end{align}
then, by using the same approach as Step~3, we can prove
\begin{align}
E_3(n)-E_3(n+1)=&-\frac{1}{\pi}q_8(n+1)-\frac{1}{\pi}\ln(n+\frac 34)+
\frac{1}{\pi}\ln(n+\frac 74)\label{Landau4-3}\\
&-\mathrm{MC}_3(n)+\mathrm{MC}_3(n+1)
+O\left(\frac{1}{n^{16}}\right),\nonumber
\end{align}
and thus, find
\begin{align*}
\kappa_3=-\frac{791896453750695892475}{691850212268234428416},
\lambda_3=\frac{79124827964452580408836456738931}
{23635681749960244849264556808320}.
\end{align*}
Similarly, using the~\emph{Mathematica} software can produce
\begin{align*}
E_3(n)-E_3(n+1)=\frac{67635788183819473349503637761583042747}
{155388490385564148205703855014831718400}\frac{1}{n^{15}}
+O\left(\frac{1}{n^{16}}\right).
\end{align*}
Finally, by Lemma 1 we have
\begin{align}
\lim_{n\rightarrow\infty}n^{14}E_3(n)=
\frac{9662255454831353335643376823083291821}{
 310776980771128296411407710029663436800\pi}:=C_3.
\end{align}
This completes the proof of Theorem 1.~\hfill\qed
\begin{thm} Let $\mathrm{MC}_2(n)$ be defined in Theorem~1. Then, for all integer $n\ge 0$, we have
\begin{align}
\frac{C_2}{(n+\frac 74)^{10}}
<G(n)-\frac{1}{\pi}\ln\left(n+\frac 34\right)-c_0-\mathrm{MC}_2(n)
<\frac{C_2}{(n+\frac 34)^{10}},\label{theorem2}
\end{align}
where $C_2=\frac{31675858150027835699}{5605686531912433139712\pi}
$.
\end{thm}

\begin{rem}
In fact, Theorem 2 implies that
$E_2(n)$ is a strictly decreasing function of $n$. In addition, it should be possible to establish many these types of inequalities  by using the same method of Theorem 2.
\end{rem}

\proof
First, it is not difficult to verify that~\eqref{theorem2}
is true for $n=0$. Hence, in the following we only need to prove that~\eqref{theorem2} holds for $n\ge 1$. For notational simplicity, we let
$D_2=\frac{158379290750139178495}{254803933268746960896}$, and
\begin{align}
E_2(n)=G(n)-\frac{1}{\pi}\ln(n+\frac 34)-c_0-\mathrm{MC}_2(n).\label{Landau42-1}
\end{align}
Then, it follows from \eqref{G-difference-def}
\begin{align}
E_2(n)-E_2(n+1)=&-\frac{1}{\pi}q(n+1)
-\frac{1}{\pi}\ln(n+\frac 34)-\mathrm{MC}_2(n)\label{Landau42-2}\\
&+\frac{1}{\pi}\ln(n+\frac 74)
+\mathrm{MC}_2(n+1).\nonumber
\end{align}
If we let
\begin{align}
U(x)=&-\frac{1}{\pi}q_8(x+1)-\frac{1}{\pi}\ln(x+\frac 34)-\mathrm{MC}_2(x)
+\frac{1}{\pi}\ln(x+\frac 74)
+\mathrm{MC}_2(x+1),\label{Landau42-3}\\
V(x)=&-\frac{1}{\pi}q_7(x+1)-\frac{1}{\pi}\ln(x+\frac 34)-\mathrm{MC}_2(x)
+\frac{1}{\pi}\ln(x+\frac 74)
+\mathrm{MC}_2(x+1),\label{Landau42-4}
\end{align}
then, combining~\eqref{Landau-1} and~\eqref{Landau42-2}-\eqref{Landau42-4} yields
\begin{align}
V(n)<E_2(n)-E_2(n+1)<U(n).\label{Landau42-5}
\end{align}
In the following, we establish the lower bound of $V(n)$ and
the upper bound of $U(n)$. First,
by using the \emph{Mathematica} software, we can obtain
\begin{align}
-U'(x)-\frac{D_2}{\pi(x+\frac 54)^{12}}=-\frac{1}{\pi}
\frac{\Psi_1(32;n)}{75937489649280 (3+4n)(5+4n)^{12}(7+4 n)\Psi_2(32;n)}<0.
\label{Landau42-6}
\end{align}
Noticing $U(+\infty)=0$, and utilizing \eqref{Landau42-6} and Lemma 4, we have
\begin{align}
U(n)&=\int_{n}^{\infty}-U'(x)dx<\int_{n}^{\infty}\frac{D_2}
{\pi(x+\frac 54)^{12}}dx
=\frac{D_2}{11\pi}\frac{1}{(n+\frac 54)^{11}}\label{Landau42-7}\\
&<\frac{D_2}{11\pi}\int_{n+\frac 34}^{n+\frac 74}\frac{1}{x^{11}}dx.\nonumber
\end{align}
Similarly, we can attain
\begin{align}
-V'(x)-\frac{D_2}{\pi(x+\frac 74)^{12}}=\frac{1}{\pi}\frac{\Psi_3(30;n)}{75937489649280(3+4 n)(5+4 n)^2(7+4n)^{12}\Psi_4(28;n)}>0.
\label{Landau42-8}
\end{align}
Therefore, integrating~\eqref{Landau42-8} with $V(+\infty)=0$ results in
\begin{align}
V(n)&=\int_{n}^{\infty}-V'(x)dx>\int_{n}^{\infty}
\frac{D_2}{\pi(x+\frac 74)^{12}}dx
=\frac{D_2}{11\pi}\frac{1}{(n+\frac 74)^{11}}\label{Landau42-9}\\
&>\frac{D_2}{11\pi}\int_{n+\frac {7}{4}}^{n+\frac {11}{4}}\frac{1}{x^{11}}dx,\nonumber
\end{align}
which, along with $E_2(\infty)=0$ and~\eqref{Landau42-5} gives us
\begin{align}
E_2(n)=&\sum_{m=n}^{\infty}\left(E_2(m)
-E_2(m+1)\right)
>\sum_{m=n}^{\infty}\frac{D_2}{11\pi}\int_{m+\frac 74}^{m+\frac {11}{4}}\frac{1}{x^{11}}dx\label{Landau42-11}\\
=&\frac{D_2}{11\pi}\int_{n+\frac 74}^{\infty}\frac{1}{x^{11}}dx=\frac{D_2}{110\pi}\frac{1}{(n+\frac 74)^{10}}.
\nonumber
\end{align}
Similarly, combining~\eqref{Landau42-7} with~\eqref{Landau42-5} yields
\begin{align}
E_2(n)=&\sum_{m=n}^{\infty}\left(E_2(m)
-E_2(m+1)\right)
<\sum_{m=n}^{\infty}\frac{D_2}{11\pi}\int_{m+\frac 34}^{m+\frac 74}\frac{1}{x^{11}}dx\label{Landau42-10}\\
=&\frac{D_2}{11\pi}\int_{n+\frac 34}^{\infty}\frac{1}{x^{11}}dx
=\frac{D_2}{110\pi}\frac{1}{(n+\frac 34)^{10}}.\nonumber
\end{align}
This finishes the proof of Theorem 2.~\hfill\qed

\section{The Lebesgue constants}
For the Lebesgue constants, we will prove the following hybrid-type finite continued fraction approximations, which has a structure  similar to that of the Landau constants.

\begin{thm}  Let the initial-correction function be given by $\mathrm{MC}_0(n)=\frac{4}{\pi^2}\ln (n+1)+c_1$, where $c_1$ is defined by~\eqref{c1-definition}. If we let the $k$-th correction
function $\mathrm{MC}_k(n)$ for $k\ge 1$ be defined by
\begin{align}
\mathrm{MC}_1(n):=&\frac{\rho_1}{(n+1)^2+\varrho_1},\\
\mathrm{MC}_k(n):=&\frac{\rho_1}{(n+1)^2+\varrho_1+}
\K_{j=2}^{k}\frac{\rho_j}{(n+1)^2+\varrho_j},
(k\ge 2).\label{LMC-definition}
\end{align}
where
\begin{align*}
\rho_1=&\frac{12-\pi^2}{18 \pi^2},\\
\varrho_1=&\frac{7(-720+60\pi^2+\pi^4)}
{600(-12+\pi^2)},\\
\rho_2=&-\frac{7515244800-1252540800\pi^2
+46937520\pi^4+65640 \pi^6
+23797\pi^8}{52920000(-12+\pi^2)^2},\\
\varrho_2=&7(-36262162944000+9065540736000\pi^2
-720128102400\pi^4+16206350400\pi^6\\
&+117169920\pi^8+288540\pi^{10}+230953\pi^{12})/
(600(-90182937600+22545734400\pi^2\\
&-1815791040\pi^4+46149840\pi^6-219924\pi^8
+23797 \pi^{10}))
\end{align*}
and the corresponding $k$-th correction error term $E_k(n)$ be defined by
\begin{align}
E_k(n):=L_{n/2}-\mathrm{MC}_0(n)-\mathrm{MC}_k(n),
\label{Lebesgue-E-definition}
\end{align}
then, for all positive integer $k$, we have
\begin{align}
&\lim_{n\rightarrow\infty}n^{4k+3}\left(E_k(n)
-E_k(n+1)\right)
=(4k+2)C_k,\\
&\lim_{n\rightarrow\infty}n^{4k+2}E_k(n)=C_k,
\label{LMC-result}
\end{align}
where
\begin{align*}
C_1=&\frac{-7515244800 + 1252540800\pi^2-46937520\pi^4-65640\pi^6-23797\pi^8}
{952560000\pi^2(-12+\pi^2)},\\
C_2=&(7633889107527073628160000
- 1908472276881768407040000\pi^2
+146687085183488661504000\pi^4\\
&-3184401328004768256000\pi^6+50811629937851059200\pi^8
-5860796365392595200\pi^{10}\\
&+73433337261096960\pi^{12}-2698623258901920\pi^{14}
-13989723377364\pi^{16}-552278517605\pi^{18})\\
&/(97592743987200\pi^2(7515244800-1252540800\pi^2
+46937520\pi^4+65640\pi^6+23797\pi^8)).
\end{align*}
\end{thm}

\proof Since the proof of Theorem 3 is very similar to that of Theorem 1, we only outline the idea of the proof here. First, we recall that
\begin{align}
E_k(n):=L_{n/2}-\mathrm{MC}_0(n)-\mathrm{MC}_k(n).
\end{align}
For every positive integer $M$, we let
\begin{align}
W_M(n):=\sum_{j=1}^{M}\frac{a_j}{(n+1)^{2j}},
\end{align}
where $a_j$ is given in Lemma 3. It is not hard to see that $a_j>0$ for odd $j=1,3,\cdots$, and $a_j<0$ for even $j=2,4,\cdots$. It follows easily from Lemma 3 and (5.6) that
\begin{align}
E_k(n)=&W_{2k+1}(n)-\mathrm{MC}_k(n)
+O\left(n^{4k+4}\right),\\
E_k(n)-E_k(n+1)=&W_{2k+1}(n)-\mathrm{MC}_k(n)-W_{2k+1}(n+1)+
\mathrm{MC}_k(n+1)
+O\left(n^{4k+4}\right).
\end{align}
Hence, it suffices for us to approximate $W_{2k+1}(n)$. Similar to the proof of Theorem~1, we expand $ W_{2k+1}(n)-\mathrm{MC}_k(n)-W_{2k+1}(n+1)+
\mathrm{MC}_k(n+1)$ into a power series in terms of $n^{-1}$, and then check \eqref{LMC-result} holds. \qed

The main purpose of this section is to prove the following theorem,
which corresponds to inequalities (1.10) in the case of Landau constants.
\begin{thm} Let $\mathrm{MC}_1(n)$ be defined in Theorem 3. Then, for all integer $n\ge 0$, we have
\begin{align}
\frac{C_1}{(n+\frac {13}{8})^6}<L_{n/2}-\frac{4}{\pi^2}\ln (n+1)-c_1-\mathrm{MC}_1(n)<\frac{C_1}{(n+\frac 58)^6},
\label{Theorem 4}
\end{align}
where $C_1=
\frac{-7515244800+1252540800\pi^2-46937520\pi^4-
65640\pi^6-23797\pi^8}{952560000\pi^2(-12+\pi^2)}>0.$
\end{thm}
\begin{rem}
In fact, Theorem 4 implies that
$E_1(n)$ is a strictly decreasing function of $n$.
\end{rem}
\proof First, since $L_{n/2}=1$ for $n=0$, one may verify that \eqref{Theorem 4}
is true for $n=0$. When $n=1$, it follows from Lemma 3 and (5.7)
that
\begin{align*}
\mathrm{MC}_0(n)+W_4(n)<L_{n/2}<\mathrm{MC}_0(n)+W_3(n),
\end{align*}
it is not difficult to verify that \eqref{Theorem 4}
is also true for $n=1$. Hence, in the following we only need to prove that \eqref{Theorem 4} holds for $n\ge 2$. By \eqref{Lebesgue-E-definition} and Lemma 3 we have
\begin{align}
&E_1(n)-E_1(n+1)\label{L-E-Difference-ub}\\
=&\left(L_{n/2}-\mathrm{MC}_0(n)\right)
-\mathrm{MC}_1(n)-
\left(L_{(n+1)/2}-\mathrm{MC}_0(n+1)\right)+\mathrm{MC}_1(n+1)
\nonumber\\
<& W_3(n)-\mathrm{MC}_1(n)-W_4(n+1)+\mathrm{MC}_1(n+1)\nonumber\\
=&\left(W_3(n)-\mathrm{MC}_1(n)-W_3(n+1)+\mathrm{MC}_1(n+1)\right)
-\frac{a_4}{(n+2)^8}.\nonumber
\end{align}
Similarly, we also have
\begin{align}
E_1(n)-E_1(n+1)>\left(W_3(n)-\mathrm{MC}_1(n)-W_3(n+1)
+\mathrm{MC}_1(n+1)\right)
+\frac{a_4}{(n+1)^8}.\label{L-E-Difference-lb}
\end{align}
For notational simplicity, we let $D_1=42C_1$. Now we define for $x\ge 1$
\begin{align}
F(x):=&W_3(x)-\mathrm{MC}_1(x)-W_3(x+1)
+\mathrm{MC}_1(x+1),\label{F(x)-Def}\\
U(x):=&\frac{D_1}{(x+\frac 54)^8},\quad V(x):=\frac{D_1}{(x+\frac 32)^8}.
\end{align}
In the following, we establish the upper bound and lower bounds of $F(x)$, respectively. By using \emph{Mathematica} software, one can check
\begin{align}
-F'(x)-U(x)=\frac{P_1(21;x)}
{19845000\pi^2(-12+\pi^2)(1+x)^7(2+x)^7(5+4x)^8\Psi_1(4;x)\Psi_2(4;x)},
\end{align}
where polynomial $P_1(21;x)$ may be expressed as
\begin{align}
P_1(21;x)=(x-1)\left(\frac{b_{-1}}{x-1}+b_0+b_1x
+\cdots+b_{20}x^{20}\right).
\end{align}
By using \emph{Mathematica} software again, it is not hard to verify that all coefficients $b_j(-1\le j\le 20)$ are positive. Thus, the inequality $P_1(21;x)>0$ holds for $x\ge 1$. Noticing that $-12+\pi^2<0$, one obtains by Lemma 4
\begin{align}
-F'(x)<&U(x),\quad x\ge 1,\\
F(n)=&\int_{n}^{+\infty}-F'(x)dx\le \int_{n}^{+\infty}U(x)dx
=\frac{D_1}{7}\frac{1}{(n+\frac 54)^7}\le
\frac{D_1}{7}\int_{n+\frac 34}^{n+\frac 74}\frac{dx}{x^{7}}.
\label{F(n)-upper bound}
\end{align}
On the other hand, we can prove by using \emph{Mathematica} software \begin{align}
-F'(x)-V(x)=\frac{P_2(20;x)}{19845000\pi^2(-12+\pi^2)(1+x)^7)(2+x)^7(3+
2x)^8 \Psi_1(4;x)\Psi_2(4;x)},
\end{align}
where
\begin{align}
P_2(20;x)=d_0+d_1x
+\cdots+d_{20}x^{20},
\end{align}
and all coefficients $d_j(0\le j\le 20)$ are negative. Thus, this yields
\begin{align}
-F'(x)>& V(x),\quad x\ge 1,\\
F(n)=&\int_{n}^{+\infty}-F'(x)dx\ge \int_{n}^{+\infty}V(x)dx
=\frac{D_1}{7}\frac{1}{(n+\frac 32)^7}\ge
\frac{D_1}{7}\int_{n+\frac 32}^{n+\frac 52}\frac{dx}{x^{7}}.
\label{F(n)-lower bound}
\end{align}
Combining \eqref{L-E-Difference-ub}, \eqref{L-E-Difference-lb},  \eqref{F(x)-Def}, \eqref{F(n)-upper bound} and  \eqref{F(n)-lower bound} gives
\begin{align}
\frac{D_1}{7}\int_{n+\frac 32}^{n+\frac 52}\frac{dx}{x^{7}}+\frac{a_4}{(n+1)^8}
<E_1(n)-E_1(n+1)<\frac{D_1}{7}\int_{n+\frac 34}^{n+\frac 74}\frac{dx}{x^{7}}-\frac{a_4}{(n+2)^8}.
\end{align}
By adding the estimates from $n$ to $\infty$ and noticing $E_1(\infty)=0$, we attain
\begin{align}
\frac{D_1}{7}\int_{n+\frac 32}^{\infty}\frac{dx}{x^{7}}
+a_4\sum_{m=n}^{\infty}\frac{1}{(m+1)^8}<E_1(n)<
\frac{D_1}{7}\int_{n+\frac 34}^{\infty}\frac{dx}{x^{7}}
-a_4\sum_{m=n}^{\infty}\frac{1}{(m+2)^8}.
\label{L-E_1(n)-double bounds}
\end{align}
By Lemma 4 again, one has
\begin{align}
\sum_{m=n}^{\infty}\frac{1}{(m+2)^8}\le\sum_{m=n}^{\infty}
\int_{m+\frac 32}^{m+\frac 52}\frac{dx}{x^8}=\frac{1}{7(n+\frac 32)^7}.
\end{align}
On the other hand, one has the following trivial estimate
\begin{align}
\sum_{m=n}^{\infty}\frac{1}{(m+1)^8}\ge \int_{n+1}^{\infty}\frac{dx}{x^8}=\frac{1}{7(n+1)^7}.
\end{align}
Substituting the above two estimates into \eqref{L-E_1(n)-double bounds} produces
\begin{align}
\frac{C_1}{(n+\frac 32)^6}+\frac{a_4}{7(n+1)^7}
<E_1(n)<
\frac{C_1}{(n+\frac 34)^6}
-\frac{a_4}{7(n+\frac 32)^7}.\label{L-E1(n)-double bounds-1}
\end{align}
By using \emph{Mathematica} software, it is not difficult to check
\begin{align}
&\frac{C_1}{(n+\frac 58)^6}-\left(\frac{C_1}{(n+\frac 34)^6}
-\frac{a_4}{7(n+\frac 32)^7}\right)\\
=&\frac{P_3(12;n)}{29767500\pi^2(-12+\pi^2)(3+2n)^7 (3+4n)^6(5+8n)^6},\nonumber
\end{align}
where
\begin{align}
P_3(12;n)=(n-1)\left(\frac{\theta_{-1}}{n-1}+\theta_0+\theta_1n+
\cdots+\theta_{11}n^{11}\right).
\end{align}
By utilizing \emph{Mathematica} software again, we observe that all coefficients
$\theta_j (-1\le j\le 11)$ are negative. Hence
\begin{align}
\frac{C_1}{(n+\frac 58)^6}-\left(\frac{C_1}{(n+\frac 34)^6}
-\frac{a_4}{7(n+\frac 32)^7}\right)\ge 0,\quad (n\ge 1).
\end{align}
Similarly, one may check
\begin{align}
&\left(\frac{C_1}{(n+\frac 32)^6}+\frac{a_4}{7(n+1)^7}\right)-
\frac{C_1}{(n+\frac {13}{8})^6}\\
=&\frac{(n-2)\left(\frac{\vartheta_{-1}}{n-2}+\vartheta_0
+\vartheta_1n+\cdots+\vartheta_{11}n^{11}
\right)}
{3810240000\pi^2(-12+\pi^2)(1+n)^7(3+2n)^6(13+8n)^6},\nonumber
\end{align}
and all coefficients $\vartheta_j(-1\le j\le 11)$ are negative.
Thus, we obtain
\begin{align}
\left(\frac{C_1}{(n+\frac 32)^6}+\frac{a_4}{7(n+1)^7}\right)-
\frac{C_1}{(n+\frac {13}{8})^6}>0,\quad (n\ge 2).
\end{align}
Finally, Theorem 4 follows from \eqref{L-E1(n)-double bounds-1}
, (5.30) and (5.32) immediately.\qed

\section{The Euler-Mascheroni constant}
The Euler constant was first introduced by Leonhard Euler (1707-1783) in 1734 as the limit of the sequence
\begin{align}
\gamma(n):=\sum_{m=1}^{n}\frac 1m -\ln n.
\end{align}
It is also known as the Euler-Mascheroni constant. There are many famous unsolved problems about the nature of this constant. See e.g. the survey papers or books of R.P. Brent and P. Zimmermann~\cite{BZ}, Dence and Dence~\cite{DD}, Havil~\cite{Hav} and Lagarias~\cite{Lag}. For example, a long-standing open problem is whether or not it is a rational number.

In fact, the sequence $\left(\gamma(n)\right)_{n\in \mathbb{N}}$ converges very slowly toward $\gamma$, like $(2n)^{-1}$. Up to now, many authors are
preoccupied to improve its rate of convergence. See e.g.~\cite{CM,DD,De,GI,GS,Lu,Lu1,Mor1,MC} and references therein.
Let $R_1(n)=\frac{a_1}{n}$ and for $k\ge 2$
\begin{align}
R_k(n):=\frac{a_1}{n+\frac{a_2 n}{n+\frac{a_3n}{n+\frac{a_4n}{\frac{\ddots}{n+a_k}}}}
}\label{Lu-CF},
\end{align}
where
$(a_1,a_2,a_4,a_6,a_8,a_{10},a_{12})=\left(
\frac{1}{2},\frac{1}{6},\frac{3}{5},
\frac{79}{126},\frac{7230}{6241},
\frac{4146631}{3833346},
\frac{306232774533}{179081182865}\right)$,
$a_{2k+1}=-a_{2k}$ for $1\le k\le 6$, and
\begin{align}
r_k(n):=\sum_{m=1}^{n}\frac 1m -\ln n-R_k(n).
\end{align}
Lu~\cite{Lu} introduced the continued fraction method to investigate
this problem, and showed
\begin{align}
\frac{1}{120(n+1)^{4}}< r_{3}(n)-\gamma<\frac{1}{120(n-1)^{4}}.
\label{Lu-inequalities}
\end{align}
Xu and You~\cite{XY} continued Lu's work
to find $a_5,\cdots, a_{13}$ with the help of \emph{Mathematica} software, and obtained
\begin{align}
\lim_{n\rightarrow\infty}n^{k+1}\left(r_k(n)-\gamma\right)=C_k',
\label{XY-results}
\end{align}
where
$(C_1',\cdots,C_{13}')=\left(-\frac{1}{12},-\frac{1}{72},\frac{1}{120},
\frac{1}{200},-\frac{79}{25200},-\frac{6241}{3175200},
\frac{241}{105840},
\frac{58081}{22018248},-\frac{262445}{91974960},\right.$

$\left. -\frac{2755095121}{892586949408},\frac{20169451}{3821257440},
\frac{406806753641401}{45071152103463200},
-\frac{71521421431}{5152068292800}\right)$.
Hence the rate of the convergence of the sequence $\left(r_k(n)\right)_{n\in \mathbb{N}}$ is $n^{-(k+1)}$. Moreover, they improved \eqref{Lu-inequalities} to
\begin{align}
&C_{10}'\frac{1}{(n+1)^{11}}<\gamma -r_{10}(n)<C_{10}'\frac{1}{n^{11}},\\
&C_{11}'\frac{1}{(n+1)^{12}}< r_{11}(n)-\gamma<C_{11}'\frac{1}{n^{12}}.
\end{align}
The purpose of this section is to further refine the works of Lu~\cite{Lu} and Xu and You~\cite{XY} by using the \emph{multiple-correction method}, and prove the following theorem.
\begin{thm}\label{th_euler}
For every positive integer $k$, let the $k$-th correction
function $\mathrm{MC}_k(n)$ be defined by
\begin{align}
\mathrm{MC}_1(n):=&\frac{a_1}{n+b_1},\\
\mathrm{MC}_k(n):=&\frac{a_1}{n+b_1+}
\K_{j=2}^{k}\frac{a_j}{n+b_j}
,(k\ge 2),
\end{align}
where
\begin{align*}
&a_1=\frac{1}{2}, \quad &b_1=\frac{1}{6},\\
&a_2=\frac{1}{36}, \quad &b_2=\frac{13}{30},\\
&a_3=\frac{9}{25}, \quad &b_3=\frac{17}{630},\\
&a_4=\frac{6241}{15876}, \quad &b_4=\frac{417941}{786366},\\
&a_5=\frac{52272900}{38950081}, \quad &b_5=-\frac{1835967509}{23923912386},\\
&a_6=\frac{17194548650161}{14694541555716}, \quad &b_6=\frac{431312596940299603}{686480136010816290},\\
&a_7=\frac{93778512198179213368089}{32070070056327569608225}, \quad &b_7=-\frac{75178865368857369613934863 }{437108607837436422694763190},\\
&a_8=\frac{14093175882028689333655328914081}
{5957702453097198927838844740836}, \quad &b_8=\frac{152838545298199920648591716358691154137
}{212256305311307139071033336233757302422}.
\end{align*}
If we let the $k$-th correction error term $E_k(n)$ be defined by
\begin{align}
E_k(n):=\sum_{m=1}^{n}\frac 1m-\ln n-\gamma -\mathrm{MC}_k(n),
\end{align}
then, for all positive integer $k$, we have
\begin{align}
&\lim_{n\rightarrow\infty}n^{2k+2}\left(E_k(n)-E_k(n+1)\right)
=(2k+1)C_k,\\
&\lim_{n\rightarrow\infty}n^{2k+1}E_k(n)=C_k,
\end{align}
where $(C_1,C_2,C_3,\cdots,C_{6})=(-\frac{1}{72},\frac{1}{200},
-\frac{6241}{3175200},\frac{58081}{22018248},
-\frac{2755095121}{892586949408},\frac{406806753641401}
{45071152103463200})$, $C_7=-\frac{5115313723510706087761}{239581189590134660611200}$ and $C_8=\frac{26329150006913625404731665769}
{241842252367746831359300280968}$.
\end{thm}
\begin{rem}
For comparison with~\eqref{Lu-CF}, Theorem~\ref{th_euler} is more convenient for us to
find $a_k$ and $b_k$, since the parameter $a_k$ and the variable $n$ in Lu's continued faction
need to be iterated more times than ours in the recursive algorithm.
\end{rem}
\begin{rem}
It is interesting to note that we have $|C_k|<1$ for all positive integers $k~(1\le k\le 8)$, and that the correction function $\mathrm{MC}_k(n)$ for all positive integers $k$ is a rational function in the form of $\frac{P_{k-1}(n)}{Q_k(n)}$ with $P_{k-1}(x),Q_k(x)\in \mathbb{Q}[x]$. Therefore, the rate of convergence for the $(E_k(n))_{k\ge 1}$ is much faster than the geometric series if we can repeat the multiple-correction infinite times. Hence, Theorem 5 predicts that it may be possible for us to find a more rapidly or
BBP-type series expansion for the Euler-Mascheroni constant in the future.
\end{rem}
\noindent{\emph{Proof of Theorem~\ref{th_euler}}}:
 This proof consists of the following steps.

\noindent{\bf (Step 1) The initial-correction.}
We choose $\mathrm{MC}_0(n)=0$, in other words, we don't need the initial-correction, and let
\begin{align}
E_0(n)=\sum_{m=1}^{n}\frac 1m -\ln n -\gamma -\mathrm{MC}_0(n)=\sum_{m=1}^{n}\frac 1m -\ln n.\label{section6-1}
\end{align}
It is not difficult to check that
\begin{align}
\lim_{n\rightarrow \infty}n^2\left(E_0(n)-E_0(n+1)\right)
=\lim_{n\rightarrow \infty}n^2\left(\ln (1+\frac 1n)-\frac{1}{n+1}\right)=\frac 12.\label{section6-2}
\end{align}
Using Lemma~1 and Noting that $E_0(\infty)=0$,  we have
\begin{align}
\lim_{n\rightarrow \infty}n E_0(n)=\frac{1}{2}=:a_1=C_0.
\label{section6-3}
\end{align}

\noindent{\bf (Step 2) The first-correction.} We let
\begin{align}
\mathrm{MC}_1(n)=\frac{a_1}{\Phi_1(1;n)}
=\frac{a_1}{n+b_1},\label{section6-4}
\end{align}
and define
\begin{align}
E_1(n)=\sum_{m=1}^{n}\frac 1m -\ln n -\gamma -\mathrm{MC}_1(n).\label{section6-5}
\end{align}
By making use of \emph{Mathemaica} software, we expand the difference $E_1(n)-E_1(n+1)$ into a power series in terms of $n^{-1}$:
\begin{align}
E_1(n)-E_1(n+1)=&\ln (1+\frac 1n)-\frac{1}{n+1}-\mathrm{MC}_1(n)+\mathrm{MC}_1(n+1)\\
=&\frac{-\frac 16+b_1}{n^3}+\frac{1-6b_1-6b_1^2}{4 n^4}+
O\left(\frac{1}{n^5}\right).\nonumber
\end{align}
By Lemma 1, the fastest sequence $E_1(n)_{n\ge 1}$ is obtained by enforcing the first
coefficient of this power series to be zero. In this case, $b_1=\frac 16$ and thus,
\begin{align}
\lim_{n\rightarrow \infty}n^4\left(E_1(n)-E_1(n+1)\right)=-\frac{1}{24}.
\label{section6-6}
\end{align}
Applying Lemma 1 again yields
\begin{align}
\lim_{n\rightarrow \infty}n^{3}E_1(n)=-\frac{1}{72}:=C_1.
\label{section6-7}
\end{align}
\noindent{\bf (Step 3) The second-correction.}
We choose
\begin{align}
\mathrm{MC}_2(n)=\frac{a_1}{n+b_1+}\frac{a_2}{n+b_2}
\end{align}\label{section6-8}
and define
\begin{align}
E_2(n)=\sum_{m=1}^{n}\frac 1m -\ln n -\gamma -\mathrm{MC}_2(n).\label{section6-9}
\end{align}
Similar to the first-correction, \emph{Mathemaica} software helps us find the power series of $E_2(n)-E_2(n+1)$ in terms of $n^{-1}$:
\begin{align}
E_2(n)-E_2(n+1)=&\ln (1+\frac 1n)-\frac{1}{n+1}-\mathrm{MC}_2(n)+\mathrm{MC}_2(n+1)
\label{section6-10}\\
=&\frac{-\frac{1}{24}+\frac{3a_2}{2}}{n^4}+\frac{\frac{17}{135}-
\frac{11a_2}{3}-2a_2 b_2}{n^5}\nonumber\\
&+\frac{-641+17820a_2-6480a_2^2+15120 a_2 b_2+6480 a_2 b_2^2}{
2592 n^6}+
O\left(\frac{1}{n^7}\right).\nonumber
\end{align}
In order to obtain the fastest convergence of the sequence from~\eqref{section6-10}, we enforce
\begin{align*}
\begin{cases}
-\frac{1}{24}+\frac{3a_2}{2}=0,\\
\frac{17}{135}-
\frac{11a_2}{3}-2a_2 b_2=0,
\end{cases}
\end{align*}
which is equivalent to
\begin{align}
a_2=\frac{1}{36} \quad \mbox{and}\quad b_2=\frac{13}{30}.
\end{align}
Therefore, we attain
\begin{align}
E_2(n)-E_2(n+1)=\frac{1}{40}\frac{1}{n^6}
+O\left(\frac{1}{n^7}\right).\label{section6-11}
\end{align}
Now by Lemma 1 again, we have
\begin{align}
\lim_{n\rightarrow \infty}n^{6}E_2(n)=\frac{1}{200}:=C_2.
\label{section6-12}
\end{align}
Since the derivations from the third-correction to the eighth-correction are very similar, here we only give the proof of the eighth-correction.

\noindent{\bf (Step 9) The eighth-correction.} We let
\begin{align}
\mathrm{MC}_8(n)=\frac{a_1}{n+b_1+}
\K_{j=2}^{8}\frac{a_j}{n+b_j},
\end{align}\label{section6-13}
and define
\begin{align}
E_8(n)=\sum_{m=1}^{n}\frac 1m -\ln n -\gamma -\mathrm{MC}_8(n).\label{section6-14}
\end{align}
Again, we resort to \emph{Mathemaica} software to expand the difference $E_8(n)-E_8(n+1)$ into a power series in terms of $n^{-1}$:
\begin{align}
&E_8(n)-E_8(n+1)\label{section6-15}\\
=&\ln (1+\frac 1n)-\frac{1}{n+1}-\mathrm{MC}_8(n)+\mathrm{MC}_8(n+1)\nonumber
\\
=&\frac{-\frac{5115313723510706087761}{15972079306008977374080}+
\frac{1220420260924203 a_8}{9014230420692640}}{n^{16}}
+O\left(\frac{1}{n^{19}}\right)
\nonumber\\
&+\frac{\phi_1}{427954568312084496528235420668724754292663281100
n^{17}}\nonumber\\
&+\frac{\phi_2}{393249450106833933773320732143901514534088027635
9961361938547200
n^{18}},\nonumber
\end{align}
where
\begin{align*}
\phi_1=&1651455193723916359597152051576300956283985319207\\
&-653628703213874127417970758864749834347105473339 b_8\\
&-61802656649700006308879054758037254288686256518 a_8 b_8\\
\phi_2=&-840200304617857764114169221449745246233826647190833037
98933256161\\
&+34897177003226720866547593123620552491030169390762369825972005484 a_8\\
&+6683314557607540897917106944093661063444345868437190779381411152 a_8b_8\\
&+603401714833886490440852089976344096105815206107927672101641232
b_8^2 a_8\\
&- 603401714833886490440852089976344096105815206107927672101641232 a_8^2.
\end{align*}
By enforcing $a$ and $b$ in~\eqref{section6-15} to satisfy the following condition:
\begin{align*}
\begin{cases}
-\frac{5115313723510706087761}{15972079306008977374080}+
\frac{1220420260924203 a_8}{9014230420692640}=0,\\
\phi_1=0,
\end{cases}
\end{align*}
i.e.,
\begin{align*}
a_8=\frac{14093175882028689333655328914081}
{5957702453097198927838844740836},
b_8=\frac{152838545298199920648591716358691154137
}{212256305311307139071033336233757302422},
\end{align*}
we obtain
\begin{align}
E_8(n)-E_8(n+1)=\frac{26329150006913625404731665769}
{14226014845161578315252957704}\frac{1}{n^{18}}
+O\left(\frac{1}{n^{19}}\right).\label{section6-16}
\end{align}
Now by Lemma 1 again, we finally attain
\begin{align}
\lim_{n\rightarrow \infty}n^{17}E_2(n)=\frac{26329150006913625404731665769}
{241842252367746831359300280968}:=C_8.
\end{align}
This completes the proof of Theorem 5.~\hfill\qed

\noindent{\bf Appendix.} For the reader's convenience, here we give some more constants in our theorems. In fact, one can use \emph{Mathematica command} ``Together" and ``Coefficient" to find
more constants.
\begin{align*}
\kappa_4=&-\frac{382149699786434954423663192287642772100258949239}{
69454986539981103777883874787703791756725862400},\\
\kappa_5=&-\frac{6793276653107510374395619138801559957929538247747
9617083158574631741
016483865590781675}{
400160623295066935399490057280725508008329690840570325893115626008
0111590006888051712},\\
\lambda_4=&\frac{3047183642643398321446537081211433153790774725
879204120678621187}{
476180753216552458418280167270798333222960626510492964456863360}.\\
\rho_3=&-25(91606669290324883537920000-     30535556430108294512640000 \pi^2 \\
&+3668717299083632345088000 \pi^4-184899901119545880576000\pi^6\\
&+3794140887258980966400\pi^8-121141186322562201600\pi^{10}\\
&+6741996412525758720\pi^{12}-105816816367920000\pi^{14}\\
&+2530746578373552\pi^{16}+7362381166104 \pi^{18}+552278517605\pi^{20})/\\
&(2561328(7515244800-1252540800\pi^2+46937520\pi^4+65640\pi^6+
23797 \pi^8)^2)\\
a_9=&\frac{38559153745620009525389781729558359566448528400
}{7562099567591782725341311886983340261624011969}, \\ b_9=&-\frac{4311810252990337765692084981855831368824641381949822699}
{16443847302827668255907904514549005300064801885045499646}\\
a_{10}=&\frac{1424408165569510820157486371128758386293642530
67475021984438416009
}{35757280329598209749962500807452853821298673049007961786630549764},
\\
b_{10}=&\frac{106368952896545249534816650756049857087954719240036868
2729425454967212
48718541
}{13163895580746317309555747820198647160355807805927403216735894475
7539476827550}.
\end{align*}

\begin{flushleft}
Xiaodong Cao\\
Department of Mathematics and Physics, \\
Beijing Institute of Petro-Chemical Technology,\\
Beijing, 102617, P. R. China \\
e-mail: caoxiaodong@bipt.edu.cn \\
\end{flushleft}

\begin{thebibliography}{99}

\bibitem{AS}M. Abramowitz and I.A. Stegun(eds.): Handbook of
Mathematical Functions with Formulas, Graphs, and Mathematical Tables. Applied Mathematics Series, vol. 55. National Bureau of Standards, Washington (1972). Ninth printing.


\bibitem{Al} H. Alzer, Inequalities for the constants of Landau
and Lebesgue, J. Comput. Appl. Math. 139(2002)215--230.



\bibitem{Ber} B.C. Berndt, Ramanujan's Notebooks, Part II, Springer-Verlag, 1999.

\bibitem{BZ}R. P. Brent and P. Zimmermann, Modern computer arithmetic. Cambridge Monographs on Applied and Computational Mathematics, 18. Cambridge University Press, Cambridge, 2011. xvi+221 pp.

\bibitem{Br} L. Brutman, A sharp estimate of the Landau constants,
J. Approx. Theory 34 (1982) 217--220.

\bibitem{CXY} X.D. Cao, H.M. Xu and X. You, Multiple-correction and faster approximation, submitted for publication(Available at: http://arxiv.org/pdf/1409.0968.pdf).


\bibitem{CC1}C.-P. Chen and J. Choi, Asymptotic expansions for the
constants of Landau and Lebesgue, Advances in mathematics, 254(2014)622--641.

\bibitem{CC2} C.-P. Chen, J. Choi, Inequalities and asymptotic expansions for the constants of Landau and Lebesgue, RGMIA Res. Rep. Collect. 17 (2014). Article 9, 11. pp. (Available at:http://rgmia.org/papers/v17/v17a09.pdf).

\bibitem{Ch1} C.-P. Chen, Approximation formulas for Landau's constants, J. Math. Anal. Appl. 387 (2012) 916--919.

\bibitem{Ch2} C.-P. Chen, Sharp bounds for the Landau constants, Ramanujan J. 31 (2013) 301--313.



\bibitem{CM} C.P. Chen, C. Mortici, New sequence converging towards the Euler-Mascheroni constant, Comput. Math. Appl. 64 (2012) 391--398.



\bibitem{CK} D. Cvijovi\'c, J. Klinowski, Inequalities for the Landau constants, Math. Slovaca, 50 (2000) 159--164.

\bibitem{CS} D. Cvijovi\'c, H.M. Srivastava, Asymptotics of the Landau constants and their relationship with hypergeometric functions, Taiwanese J. Math. 13 (2009) 855--870.




\bibitem{DD} T.P. Dence, J.B. Dence, A survey of Euler's constant, Math. Mag. 82 (2009) 255--265.

\bibitem{De} D.W. DeTemple, A quicker convergence to Euler's constant, Amer. Math. Monthly 100 (5) (1993) 468--470.

\bibitem{EF} A. Eisinberg, G. Franz\`e and N. Salerno, Asymptotic expansion and estimate of the Landau constant, Approx. Theory Appl. (N.S.) 17 (2001) 58--64.

\bibitem{Fal} L.P. Falaleev, Inequalities for the Landau constants, Sib. Math. J. 32 (1991) 896--897.

\bibitem{Fe} L. Fej\'er, Lebesguesche Konstanten und divergente
Fourierreihen, J. Reine Angew. Math. 138 (1910) 22--53.

\bibitem{Gal} P.V. Galkin, Estimates for the Lebesgue constants,
Proc. Steklov Inst. Math. 109 (1971) 1--4.

\bibitem{GI} I. Gavrea and M. Ivan,  Optimal rate of convergence for sequences of a prescribed form. J. Math. Anal. Appl. 402 (2013), no. 1, 35--43.

\bibitem{GS} X. Gourdon and P. Sebah, Collection of formulae for the Euler constant. http://numbers.computation.
free.fr/Constants/Gamma/gammaFormulas.pdf




\bibitem{Gra} H. Granath, On inequalities and asymptotic expansions for the Landau constants, J. Math. Anal. Appl. 386 (2012) 738--743.

\bibitem{Gro} T.H. Gronwall, \"uber die Lebesgueschen Konstanten
bei den Fourierschen Reihen, Math. Ann. 72 (1912) 244--261.

\bibitem{Ha} G.H. Hardy, Note on Lebesgue's constants
in the theory of Fourier series, J. Lond. Math. Soc. 17 (1942) 4--13.


\bibitem{Hav}J. Havil, Gamma: Exploring Euler's Constant, Princeton University Press, Princeton, NJ, 2003.


\bibitem{Lag} Jeffrey C. Lagarias,  Euler's constant: Euler's work and modern developments. Bull. Amer. Math. Soc. (N.S.) 50 (2013), no. 4, 527--628.

\bibitem{La} E. Landau, Absch\"atzung der Koeffzientensumme einer Potenzreihe, Arch. Math. Phys. 21 (42-50) (1913) 250--255.

\bibitem{Le} H. Lebesgue, Le\c{c}ons sur les s\'eries
Trigonom\'etriques, Gauthier-Villars, Paris, 1906.


\bibitem{Lu} Dawei Lu, A new quiker sequence convergent to Euler's
constant, J. Number Theory, 136(2014), 320--329.

\bibitem{Lu1} Dawei Lu , Some quicker classes of sequences convergent to Euler's constant. Appl. Math. Comput. 232 (2014), 172--177.

\bibitem{Mor1} C. Mortici, On new sequences converging towards the Euler-Mascheroni constant, Comput. Math. Appl. 59 (8) (2010) 2610--2614.

\bibitem{Mor2} C. Mortici, Product approximations via asymptotic integration, Amer. Math. Monthly, 117 (5) (2010) 434--441.

\bibitem{Mor3} C. Mortici, New approximations of the gamma function in terms of the digamma function,
    Applied Mathematics Letters, 23 (2010) 97--100.

\bibitem{Mor4} C. Mortici, Sharp bounds of the Landau constants, Math. Comp. 80 (2011) 1011--1018.



\bibitem{MC} C. Mortici and  C.-P. Chen,
On the harmonic number expansion by Ramanujan, J. Inequal. Appl. 2013, 2013:222, 10 pp.


\bibitem{Ne1} G. Nemes, A. Nemes, A note on the Landau constants, Appl. Math. Comput. 217 (2011) 8543--8546.

\bibitem{Ne2} G. Nemes, Proofs of two conjectures on the Landau constants, J. Math. Anal. Appl. 388 (2012) 838--844.


\bibitem{Po} E.C. Popa, Note of the constants of Landau, Gen. Math. 18 (2010) 113--117.

\bibitem{PS} E.C. Popa and N.-A. Secelean, On some inequality for the Landau constants, Taiwanese J. Math. 15 (2011) 1457--1462.


\bibitem{Sz} G. Szeg\"o, \"Uber die Lebesgueschen Konstanten
bei den Fourierschen Reihen, Math. Z. 9 (1921) 163--166.

\bibitem{Wal} J. Wallis, Arithmetica Infinitorum, Oxford, England, 1656; Facsimile of relevant pages available in: J.A. Stedall, Catching Proteus: The collaborations of
Wallis and Brouncker. I. Squaring the circle, Notes and Records Roy. Soc. London 54 (3) (2000) 293--316.

\bibitem{Wat} G.N. Watson, The constants of Landau and Lebesgue, Quart. J. Math. Oxford Ser. 1 (1930) 310--318.

\bibitem{Wo} R. Wong, Asymptotic Approximations of Integrals, Classics in Applied Mathematics, 34. Society for Industrial and Applied Mathematics (SIAM), Philadelphia, PA, 2001. xviii+543 pp. ISBN: 0-89871-497-4


\bibitem{XY} Hongmin Xu and Xu You, Continued fraction inequalities for Euler-Mascheroni constant,  J. Inequal. Appl. 2014, 2014:343 ,11 pp.


\bibitem{Ya} Shijun Yang, On an open problem of Chen and Mortici concerning the Euler-Mascheroni constant, J. Math. Anal. Appl. 396 (2012) 689--693.

\bibitem{Zhao} D. Zhao, Some sharp estimates of the constants of Landau and Lebesque, J. Math. Anal. Appl. 349 (2009) 68--73.
\end{thebibliography}
\end{document}